\documentclass[11pt,a4paper,twoside]{amsart}
\usepackage{amsmath}
\usepackage{amssymb}
\usepackage{latexsym}

\newcommand{\co}{\colon\thinspace}    %  Colon with correct spacing for maps.
\newcommand{\fnote}[1]{\footnote{\small sharp1}}
\newcommand{\N}{{\mathbb N}}
\newcommand{\Z}{{\mathbb Z}}
\newcommand{\R}{{\mathbb R}}
\newcommand{\Q}{{\mathbb Q}}

\newcommand{\T}{{\mathbb T}}

\newcommand{\inter}{\mbox{Int}}

\newtheorem{theorem}{Theorem}[section]
\newtheorem{proposition}[theorem]{Proposition}
\newtheorem{corollary}[theorem]{Corollary}

\newtheorem{lemma}[theorem]{Lemma}
\newtheorem{remark}[theorem]{Remark} 
 
\newtheorem{problem}[theorem]{Problem} 

\title{Differentiability of Mather's beta function in low dimensions}
\author{Daniel Massart}
\date{\today}
\begin{document}

\begin{abstract}
Let $L$ be a time-periodic Tonelli Lagrangian on a two-torus. Then the $\beta$-function of $L$ is differentiable in at least $k$ directions at any $k$-irrational homology class, for $k=0,1,2$. 
\end{abstract}
\maketitle
%%%%%%%%%%%%
\section{Introduction}
%%%%%%%%%%%
This paper addresses the problem of the differentiability of Mather's $\beta$-function for time-periodic Lagrangian systems. 
The setting is the dynamics of time-periodic Lagrangian systems as introduced by Mather in \cite{Mather91}. In the sequel,  $M$ is a closed, connected manifold. A Tonelli Lagrangian on $M$ is a $C^{2}$ function on $TM\times \T$, where $\T$ is the unit circle, satisfying the following conditions : 
\begin{enumerate}
\item
for every $(x,t) \in M \times \T$, the function $v \mapsto L(x,v,t)$ 
is superlinear
\item
for every $(x,v,t) \in TM \times \T$, 
the bilinear form $\partial^2 L(x,v,t) / \partial v^2$ is positive definite 
\item
the local flow $\Phi_t$  defined on $TM\times \T$  
by the Euler-Lagrange equation for extremals of the action of curves 
is complete.
\end{enumerate}

  A good  example to keep in mind is the sum of a Riemann  metric, viewed as a quadratic function on $TM$, and a time-periodic potential (a function on $M\times \T$).   See \cite{Fathi}  for more background and references.

 Define $\mathcal{M}_{inv}$
to be the set of $\Phi_t$-invariant, compactly supported, Borel probability measures on $TM\times \T$.
Mather showed that the function (called action of the Lagrangian on measures)
	\[
	\begin{array}{rcl}
\mathcal{M}_{inv} & \longrightarrow & \R \\
\mu & \longmapsto & \int_{TM\times \T}	L d\mu
\end{array}
\]
is well defined and has a minimum.  A measure achieving this minimum is called $L$-minimizing.

When $M=\T$, by Mather's Graph Theorem (\cite{Mather91}) an invariant measure has a rotation number just like an invariant measure of a circle homeomorphism. For other manifolds Mather proposed in \cite{Mather91} the following generalization. 
First he observed that if $\omega$ is a closed one-form on $M$ and $\mu \in \mathcal{M}_{inv}$ then the integral $\int_{TM\times \T}	\omega d\mu$ is well defined, and only depends on the cohomology class of $\omega$. By duality this 
defines a homology class $\left[\mu\right]$. This $\left[\mu\right]$ is the unique  $h \in H_1 (M,\R)$ such that  
	\[
\langle h,\left[\omega \right]\rangle = \int_{TM\times \T}	\omega d\mu 
\]
for any closed one-form $\omega$ on $M$. As Mather showed in \cite{Mather91}, for any $h \in H_1 (M,\R)$, the set 
	\[ \mathcal{M}_{h,inv}:= \left\{\mu \in \mathcal{M}_{inv} \co \left[\mu\right]=h\right\}
\]
is not empty. Again the action of the Lagrangian on this smaller set of measures has a minimum, which is a function of $h$, 
called the $\beta$-function of the system :
\[
	\begin{array}{rcl}
\beta \co H_1 (M,\R) & \longrightarrow & \R \\
h & \longmapsto & 
\min \left\{\int_{TM\times \T}	Ld\mu \co \left[\mu\right]=h\right\}.
\end{array}
\] 
A measure $\mu$ such that $\int L d\mu = \beta (\left[\mu\right])$ is called $(L,h)$-minimizing. 

There is a dual construction to that of $\beta$ : if $\omega$ is a closed one-form on  $M$, 
then $L-\omega$ is a Lagrangian to which Mather's theory applies, 
and furthermore $L-\omega$ has the same Euler-Lagrange flow as $L$. 
The minimum over $\mathcal{M}_{inv}$ of $\int (L-\omega)d\mu$ is actually a function of the cohomology class of $\omega$, 
the opposite of which is called the $\alpha$-function of the system :
\[
	\begin{array}{rcl}

\alpha \co H^1 (M,\R) & \longrightarrow & \R \\
c & \longmapsto & 
-\min \left\{\int_{TM\times \T}	(L-\omega)d\mu \co \mu \in  \mathcal{M}_{inv},\  \left[\omega\right]=c\right\}.
\end{array}
\]
An  $(L-\omega)$-minimizing measure is also called $(L,\omega)$-minimizing 
or $(L,c)$-minimizing if $c$ is the cohomology of $\omega$.

Mather proved that $\alpha$ and $\beta$ are convex, superlinear, and Fenchel dual of one another, that is,
\begin{eqnarray*}
\forall h \in H_1 (M,\R),\ \beta (h) &=& \sup_{c \in H^1 (M,\R)}\left( \langle c,h \rangle -\alpha (c) \right)\\
\forall c \in H^1 (M,\R),\ \alpha (c) &=& \sup_{h \in H_1 (M,\R)}\left( \langle c,h \rangle -\beta (h) \right).
\end{eqnarray*}
In particular $\min \alpha = -\beta(0)$. 
The main geometric features of a convex function are its smoothness and strict convexity, or lack thereof. 
In the present setting they turn out to have interesting dynamical meanings as well. The prototype of all theorems in the subject is 
\begin{theorem}[\cite{Mather90, Bangert94}]\label{Mather}
If $M=\T$ then $\beta$ is differentiable at every irrational homology class. 
It is differentiable at a rational homology class if and only if periodic orbits in this class fill up $\T$. 
\end{theorem}
Since $H_1 (\T,\R)=\R$ the word rational is self-explanatory. Our purpose in this paper is to extend Theorem \ref{Mather} to the next degree of freedom, that is, $M=\T^2$, so  we need a bit of terminology. 
The torsion-free part of $H_1 (M,\Z)$ embeds as a lattice $\Gamma$ in $H_1 (M,\R)$. A class 
$h \in H_1 (M,\R)$ is called integer if it lies in $\Gamma$, and rational if $nh \in \Gamma$ for some $n \in \Z$. 
A subspace of $H_1 (M,\R)$ is called integer if it is generated by integer classes.

We  need to give a quantitative meaning to the irrationality of a homology class. 
The quotient $H_1 (M,\R)/\Gamma$ is a torus $\T^{b}$, where $b$ is the first Betti number of $M$. 
For $h$ in $H_1 (M,\R)$, the image of $\Z h$ in $\T^{b}$ is a subgroup of $\T^{b}$, 
hence its closure $\mathcal{T}(h)$ is a finite union of tori of equal dimension. This dimension is called the irrationality 
$I_{\Z}(h)$ of $h$. 
It is zero if $h$ is rational. We say a class $h$ is completely irrational 
if its irrationality is maximal, i.e. equals $b$. 
In the same way, if $v$ is a vector of $\R^n$, we call irrationality of $v$ 
the dimension of the image of $\Z v$ in $\R^n / \Z^n$. 
Note that the irrationality of $h$ equals that of $nh$ for $n\in \Z$, 
$n \neq 0$ since the quotient of $\mathcal{T}(h)$ 
by $\mathcal{T}(nh)$ is a group of cardinality $n$. 

A convex function has a tangent cone at every point. We say that $\beta$ is differentiable in $k$ directions at $h$ if the tangent cone to $\beta$ at $h$ 
contains a linear space of dimension $k$. 
We are thus led to ask \textit{whether $\beta$ is always differentiable in $k$ directions at a $k$-irrational homology class}.  This will henceforth be referred to as the Differentiability Problem.
Mather conjectures the answer is yes for $C^{\infty}$ Lagrangians. 
The answer to the Differentiability Problem is yes   for all $C^{2}$ Lagrangians when $M= \T$ by Theorem \ref{Mather}. It cannot be yes in general by \cite{BIK}.   Our main result here is that the answer is yes for all $C^{2}$  Lagrangians when $M=\T^2$ : 
 \begin{theorem}\label{principal}
 Let\begin{itemize}
  \item $L \co T\T^2 \times \T  \longrightarrow \R$ be a Tonelli Lagrangian
  \item $h_0$ be a $k$-irrational homology class in $H_1 (M,\R)$, with $k=0,1,2$.
\end{itemize} 
 Then $\beta_L$ is differentiable at $h_0$ in at least $k$ directions.
 \end{theorem} 
\begin{remark} In contrast with Mather's theorem \ref{Mather}, in general it is unclear what can be said when $\beta_L$ is differentiable at a $k$-irrational homology class $h$ in $p$ directions, with $p>k$. However, when the Lagrangian is autonomous, and the homology class is non-singular (which means that its Aubry set contains no fixed point), Theorem 2 of \cite{MS} says that $\T^2$ is foliated by closed extremals with homology $h$.
\end{remark}
%%%%%%%%%%%%%%
\subsection{Historical remarks,  and open questions}
%%%%%%%%%%
In \cite{ijm} an affirmative answer is claimed to the Differentiability problem when $L$ is an autonomous Lagrangian on a manifold of dimension two. However the proof is full of gaps, and only works when the Lagrangian $L$ is a Riemannian metric of negative curvature. The case when $M=\T^2$ and $L$ is autonomous is now a particular case of Theorem \ref{principal}. In \cite{vmb2} a particular case of the differentiability problem is treated : that is, when the non-differentiability is maximal, i.e. $\beta$ is differentiable in exactly zero direction at some homology class $h$. We then say that $\beta$ has a vertex at $h$. In \cite{vmb2} we prove that if $L$ is a Tonelli Lagrangian on a manifold of dimension two, and $\beta_L$ has a vertex at some homology class $h$, then $h$ must be rational, i.e. $0$-irrational. What we do here is to tackle the intermediate degree of non-differentiability.

The logical next step is to consider surfaces of higher genus, and then manifolds of higher dimensions. Then we encounter the following question, which we believe to be interesting in its own right : 
\begin{problem}\label{probleme}
Let $P$ be a completely irrational hyperplane of $\R^n$ (that is, $P$ does not contain any integer point other than the origin), and let $\gamma \co \R \longrightarrow P$ be a continuous map such that $\gamma (t)/t$ has a limit $l \in P \setminus \{0\}$. Is it true that the closure of the image of $\gamma(\R)$ in the torus $\T^n$ has Hausdorff dimension at least two ?
\end{problem}
Here is an even simpler version of this problem : 
\begin{problem}
Let $a$ and $b$ be real numbers which are independent over $\Q$, and let $u_n$ be a sequence of real numbers such that for any $n \in \N$, we have either $u_{n+1} = u_n +a$, or $u_{n+1} = u_n +b$. Let $S$ be the closure in the circle $\T$, of the set of all values of the sequence $u_n$, modulo one. Is it true that $S$ always has positive Lebesgue measure ? 
\end{problem}
The interested reader may want to have a look at \cite{hab}, Annexe B, for a tentative discussion of these problems. 
 
\section{Aubry sets and faces of $\alpha$}
%%%%%%%%%
\subsection{Aubry sets}
We refer the reader to \cite{Fathi} for the definition of the Peierls barrier and the Aubry set associated with the Lagrangian $L$. All we need to know is that\begin{itemize}
  \item the Peierls barrier is a Lipschitz map from $(M \times \T)^2 $ to $\R$, which we denote $h_L \left( (x,t), (y,s) \right)$
  \item the Aubry set is a compact subset of $TM \times \T$, which is invariant under the Euler-Lagrange flow of $L$, which we denote $\mathcal{A}(L)$.
\end{itemize}

If $\omega$ is a closed 1-form on $M$, then $L-\omega$ is a Tonelli Lagrangian, and its  Aubry set only depends on the cohomology class $c$ of $\omega$. We denote it $\mathcal{A}(L,c)$.

Theorem \ref{principal} will come as a corollary of our next result, which links the differentiability of the $\beta$-function with the topology of the complement of the Aubry set $\mathcal{A}(L,c)$. 

\subsection{Definition of $E_c$}
We call 
\begin{itemize}
  \item supporting subspace to the graph of $\alpha$, any affine subspace of \newline
  $H^1(M,\R) \times \R$ that meets the graph of $\alpha$  but not the open epigraph 
    \[\left\{(x,t)\in H^1(M,\R) \times \R \co t> \alpha(x),        \right\}\]  
  \item  flat of   $\alpha$,  the intersection of the graph of $\alpha$ with a supporting subspace.
  \end{itemize}
Note that flats of   $\alpha$ are convex since $\alpha$ is convex, so we may speak of their relative interiors. 

Throughout this paper we view $TM \times \T$ as embedded into $T\left( M \times \T \right)$ by the map $(x,v,t) \longmapsto (x,v,t,1)$. This allows us to compare the Aubry set, which is a subset of $TM \times \T$, with the support of 1-forms on $T\left( M \times \T \right)$. 
Now if $c \in  H^1(M,\R)$, we define the following subsets of  $H^1(M,\R)\times \R$ :
\begin{itemize}
  \item $\tilde{ F}_c (L):= \{ \left(c',\alpha(c')\right) \co \mathcal{A}(L,c) \subset \mathcal{A}(L,c') \}$
  \item $F_c (L)$ is  the canonical projection of $\tilde{ F}_c (L)$ to $H^1(M,\R)$
     \item $\tilde{V}_c (L):= \{ \lambda \left( c'-c,\alpha(c)-\alpha(c')\right) \co \lambda \in \R,\  c' \in F_c \} $
  \item $\tilde{E}_c (L)$ as the set of cohomology classes of closed one-forms on $M\times \T$ which are supported outside  $\mathcal{A}(L,c)$
  \item $V_c (L), E_c (L)$ are the canonical projections of $\tilde{V}_c(L), \tilde{E}_c(L)$, respectively, to $H^1(M,\R)$.
  \end{itemize}
  We shall abbreviate the notations to $\tilde{ E}_c, \tilde{ V}_c, \tilde{ F}_c$ when there is no ambiguity on the Lagrangian.
 It can be seen from \cite{ijm}, \cite{soussol}, that $\tilde{ F}_c$ is the maximal flat of $\alpha$ containing $c$ in its relative interior. Moreover, by  \cite{vmb2}, Proposition 21 (the autonomous case of which is \cite{ijm}, Proposition 6), for any $c'$ such that  $(c',\alpha(c'))$ lies in the relative interior of $ \tilde{ F}_c$, we have $\mathcal{A}(L,c) = \mathcal{A}(L,c') $. The set $\tilde{ V}_c$ is the underlying vector space to the affine subspace of $H^1(M,\R)\times \R$ generated by $\tilde{ F}_c$.

  We proved in \cite{ijm} that $E_c \subset V_c$ for any autonomous Tonelli Lagrangian on a closed manifold $M$, and any cohomology class $c \in H^1(M,\R)$. The time-periodic case is treated in \cite{soussol}. In other words, if you have a closed one-form $\omega$ supported away from your Aubry set, you may add a small multiple of $\omega$ to your Lagrangian without changing the Aubry set. 
  
Here we prove the opposite inclusion when $M= T^2$ :
\begin{theorem}  \label{tore de dim 2}
Let
\begin{itemize}
  \item $L \co T\T^2 \times \T  \longrightarrow \R$ be a Tonelli Lagrangian
  \item $c$ be any cohomology class in $H^1 (M,\R)$.
\end{itemize} 
 Then $ E_{c}= V_{c}$, and $ V_{c}$ is an integer subspace of $H^1 (M,\R)$.
 \end{theorem}

%%%%%%%%%%%%%
\subsection{Proof of Theorem \ref{principal} assuming Theorem \ref{tore de dim 2} }
  Proposition A.3 of \cite{vmb2} reads
  \begin{proposition}\label{vmb2, A.3}
Let $L$ be a Tonelli Lagrangian on a closed manifold $M$.
Assume that for every cohomology class $c$, $\tilde{V}_c$ is an integer subspace of 
$H^1 (M\times \T,\R)$. Let $h$ be a $k$-irrational homology class. 
Then $\beta_L$ is differentiable at $h$ in at least $k$ directions.
\end{proposition}
Now Corollary 10 of \cite{vmb2} says that if 
 $L$ is a Tonelli Lagrangian on a closed manifold of dimension two,  if $V_c$ is integer, then so is $\tilde{V}_{c}$.
Since Theorem \ref{tore de dim 2} says that $V_{c}$ is an integer subspace of $H^1 (M,\R)$, Theorem \ref{principal} follows.
\qed

%%%%%
\section{Proof of Theorem \ref{tore de dim 2}}\label{preuve}
%%%%%%
Replacing $L$ with $L-\omega$, where $\omega$ is a closed 1-form with cohomology $c$, we may assume $c= 0$. For simplicity we denote $\mathcal{A} :=\mathcal{A}(L,0)$. There are two cases : 
\begin{enumerate}
  \item either  for any $c \in V_0$ there exists $\lambda \in \R^*$ such that $\lambda c \in H^1(\T^2, \Z)$
  \item or there exists $c \in V_0$ such that for all  $\lambda \in \R^*$,  $\lambda c \not\in H^1(\T^2, \Z)$.
\end{enumerate}
In the first case we observe that   $V_0$ has to be a one-dimensional subspace generated by some element of $H^1(M,\Z)$, that is, $V_0$ is an integer subspace of $H^1(M,\R)$.  Then we use  Corollary 16 of \cite{vmb2}  : 
\begin{corollary}
If the dimension of $M$ is two and $V_0$ contains an integer point $c$, then $c \in E_0$. In particular, if the dimension of $M$ is two and $V_0$ is integer, then $E_0=V_0$.
\end{corollary}
So the theorem holds in the first case. 

In the second case, we show that $E_0 = H^1(\T^2, \R)$. Since $E_0 \subset V_0$,  $E_0 = H^1(\T^2, \R)$ entails $V_0 =E_0 = H^1(\T^2, \R)$, and since $H^1(\T^2, \R)$ is an integer subspace of itself, $V_0$ is also integer. 

So in the remainder of the proof we take $c \in V_0$ and  assume that $c$ is  2-irrational, that is, for any non-zero $\lambda \in \R$, $\lambda c \not\in H^1(\T^2, \R)$, and we  prove that $E_0 = H^1(\T^2, \R)$.

Let $p \co \R^2 \longrightarrow \T^2$ be the universal cover of $\T^2$. For brevity we denote  $\overline{\mathcal{A}} := p^{-1}(\mathcal{A}) $.
Coordinates $(x,y)$ are meant with respect to the canonical basis of $\R^2$. 
Let $(\lambda, \mu)$ be the coordinates of $c$ in the basis $(\left[dx\right],\left[dy\right])$ of $H^1 (\T^2, \R)$. Since $c$ is 2-irrational,  $\lambda$ and $\mu$ are independant over $\Q$, in particular neither of them is zero. 
Define 
\begin{eqnarray*}
\omega & := &  \lambda dx +\mu dy \\
u_0 (x,y) &:=&  h_L \left( \left((x,y),0\right), (0,0) \right) \\
u_1 (x,y) &:=&  h_{L-\omega} \left( \left((x,y),0\right), (0,0) \right).
\end{eqnarray*}
Then $\omega$ is a smooth 1-form on $\T^2$ with cohomology $c$ and $u_0, u_1$ are Lipschitz functions on $\T^2$.
For simplicity we shall use the same notation for the lifts of $u_0$ and $u_1$ to $\R^2$. 
Consider the maps 
\begin{eqnarray*}
\overline{\phi} \co \R^2 & \longrightarrow &  \R^2 \\
(x,y) & \longmapsto & \left( x +  \frac{1}{\lambda} (u_1(x)-u_1(0) -u_0(x) + u_0 (0)), y \right) \\
\overline{\phi}_c \co \R^2 & \longrightarrow &  \R \\
(x,y) & \longmapsto &  \lambda x + \mu y +  u_1(x)-u_1(0) -u_0(x) + u_0 (0)\\
\pi  \co \R^2 & \longrightarrow &  \R \\
(x,y) & \longmapsto &  \lambda x + \mu y.
\end{eqnarray*}

Here are a few observations about the maps $\overline{\phi}$ and $\overline{\phi}_c$ : 
\begin{itemize}
  \item since the homology classes $\left[dx\right]$ and $\left[dy\right]$ are integer, the map $\overline{\phi}$ quotients to a map 
  $\phi \co \T^2 \longrightarrow \T^2$
  \item since the maps $u_0$ and $u_1$ are $\Z^2$-periodic, $\overline{\phi}$ is $\Z^2$-equivariant, that is, 
$$
\forall x, y \in \R, \  \forall n,m \in \Z,\   \overline{\phi} (x+n, y +m) = \overline{\phi} (x, y ) + (n,m).
$$
 As a consequence , $\overline{\phi}$ is the identity on $\Z^2$, hence $\phi$   induces the identity of $H_1 (\T^2, \R)$
  \item by \cite{vmb2}, Proposition 6, the restriction of $\overline{\phi}_c$ to $\overline{\mathcal{A}} $  satisfies a   H\H{o}lder condition of order two
  \item $\overline{\phi}_c = \pi \circ \overline{\phi}$.
\end{itemize}
The reason why H\H{o}lder  estimates on   $\overline{\phi}_c$ are interesting is  Lemma A.1 of \cite{FFR}, which originates in \cite{Ferry} :
\begin{lemma}[Ferry] \label{lemme de Ferry}
Let $\Phi$  be a map from a subset $E$ of $\R^n$ to a metric space $(X,d)$. Suppose there exist $p >1$ and $C$ such that 
	\[\forall x,y \in E,\  d\left(\Phi(x),\Phi(y)\right)\leq C \|x-y\|^p.
\]
Then the $n/p$-dimensional Hausdorff measure of $\left(\Phi (E),d\right)$ is zero.
\end{lemma}
Therefore the Hausdorff 1-dimensional measure of $\overline{\phi}_c \left( \overline{\mathcal{A}}  \right) $ is zero, so the restriction of $\overline{\phi}_c$ to $\overline{\mathcal{A}} $ is not onto. Since $\overline{\phi}_c = \pi \circ \overline{\phi}$, and the kernel of $\pi$ is the straight line $D_0$ defined by the equation $\lambda x + \mu y = 0$, it follows that there exists a straight line $D$ parallel to $D_0$, such that 
$$
\overline{\phi} \left( \overline{\mathcal{A}} \right)  \cap D = \emptyset.
$$
Now $\overline{\mathcal{A}} $ is invariant by integer translations in $\R^2$, and $\overline{\phi}$ is $\Z^2$-equivariant, so $\overline{\phi} \left( \overline{\mathcal{A}} \right) $ is also invariant by integer translations. Thus, denoting by $ \tau_{v }$ the translation by the vector $v \in \R^2$, 
$$
\forall  n,m \in \Z,\ \overline{\phi} \left( \overline{\mathcal{A}}  \right)  \cap \tau_{(n,m) }(D)= \emptyset.
$$
Recall that $c$ is 2-irrational, so $\lambda$ and $\mu$ are independant over $\Q$, hence the integer translates of $D$ are dense in $\R^2$. We identify   $\Z^2$ with $H_1 (\T^2, \Z)$ in such a way that $(1,0),(0,1)$ is the dual basis to $\left[dx\right],\left[dy\right]$. Denote by $d$ the Euclidean distance in $\R^2$. 
\begin{lemma}\label{D, epsilon}
There exists $\epsilon _0 >0$ such that for any $0< \epsilon \leq \epsilon_0$, for any $h \in  H_1 (\T^2, \Z) \cong \Z^2$ such that $d (D, D+h) \leq \epsilon$, there exists a perpendicular segment from $D$ to $D+h$ which does not meet $\overline{\phi} \left( \overline{\mathcal{A}} \right) $.
\end{lemma}
\proof
By contradiction. Assume that for any $\epsilon _0 >0$ there exists  $0< \epsilon \leq \epsilon_0$, and $h \in H_1(\T^2, \Z)$, such that 
$d (D, D+h) \leq \epsilon$, and  any perpendicular segment from $D$ to $D+h$ contains a point of $\overline{\phi} \left( \overline{\mathcal{A}}  \right) $.
Then any point of $D$ lies within distance at most $\epsilon$ of some point of $\overline{\phi} \left( \overline{\mathcal{A}}  \right) $. Therefore, since the projection $p$ does not increase distances,  any point of $p(D)$ in $\T^2$ lies within distance at most $\epsilon$ of some point of 
$\phi\left(\mathcal{A} \right) $. But $p(D)$ is dense in $\T^2$ because the integer translates of $D$ are dense in $\R^2$. Hence any point of $\T^2$ lies within distance at most $2\epsilon$ of some point of $\phi\left(\mathcal{A} \right) $. Since $\epsilon$ is arbitrarily small, it follows that 
$\phi\left(\mathcal{A} \right) $ is dense in $\T^2$. But $\mathcal{A}$ is compact, so $\phi\left(\mathcal{A} \right) $ is compact, therefore 
$\phi\left(\mathcal{A} \right) = \T^2$. Thus  $\overline{\phi} \left( \overline{\mathcal{A}} \right) = \R^2$, whence $\overline{\phi}_c \left( \overline{\mathcal{A}}  \right) = \R$, which contradicts Ferry's lemma.
\qed

By  Lemma \ref{D, epsilon}, for any $0< \epsilon \leq \epsilon_0$, and for any $h \in H_1(\T^2, \Z)$, such that 
$d (D, D+h) \leq \epsilon$, we can find  a piecewise smooth arc in $R^2$ which connects some point $x \in D$ with $x+h$ without meeting 
$\overline{\phi} \left( \overline{\mathcal{A}} \right) $. The projection of this arc to $\T^2$ is a closed curve $\gamma$ with homology $h$ which does not meet $\phi\left(\mathcal{A} \right) $. Since $\phi$ induces the identity of $H_1 (\T^2, \Z)$, $\phi^{-1} (\gamma)$ is a closed curve with homology $h$ which does not meet $\mathcal{A}$. Since $\mathcal{A}$ and $\phi^{-1} (\gamma)$ are compact, there exists a neighborhood $U$ of $\phi^{-1} (\gamma)$ which does not meet $\mathcal{A}$. By a classical construction (see for instance \cite{FK})  there exists a smooth closed 1-form $\eta_h$ supported in $U$, whose cohomology class  $c_h$ is defined by 
$$
\begin{array}{ccc}
H_1 (T^2, \R) & \longrightarrow & \R \\
k & \longmapsto & \inter (h, k)
\end{array}
$$
where $\inter$ denotes the symplectic form on $H_1 (T^2, \R) $ induced by the algebraic intersection of curves in $\T^2$. Since the 1-form $\eta_h$ is supported outside $\mathcal{A}$,  the cohomology class $c_h$ lies in $E_0$. 

Let us denote by $S_{\epsilon}$ the set of  homology classes $h \in H_1(\T^2, \Z)$, such that 
$d (D, D+h) \leq \epsilon$, and we can find  a piecewise smooth arc in $R^2$ which connects some point $x \in D$ with $x+h$ without meeting 
$\overline{\phi} \left( \overline{\mathcal{A}} ) \right) $. Denote by $S'_{\epsilon}$ the set of cohomology classes $c_h$, for 
 $h \in S_{\epsilon}$. Assume for a moment that  $S'_{\epsilon}$ generates $H^1(\T^2, \R)$. Recall that  we have seen $S'_{\epsilon} \subset  E_0$ when $\epsilon \leq \epsilon_0$. Therefore $E_0 = H^1(\T^2, \Z)$, which proves Theorem \ref{tore de dim 2}. So what we have to do now is to prove that $S'_{\epsilon}$ generates $H ^1 (M,\R)$.

Observe that the elements of $S'_{\epsilon}$ lie in the integer lattice $H^1(\T^2, \Z)$, so either they generate $H^1(\T^2, \R)$, or they lie in a subspace of dimension one generated by an element of $H^1(\T^2, \Z)$. In that case they all have the same kernel (as linear forms on $H_1(\T^2, \R)$), and that kernel is a straight line of $\R^2$ generated by an integer point. Observe that the kernel of the cohomology class $c_h$ is the straight line generated by $h$. So it is equivalent to say that the cohomology classes in $S'_{\epsilon}$ have the same kernel, and to say that there exists some 
$h_0$ in $H^1(\T^2, \Z)$ such that for any $h \in S_{\epsilon}$, there exists $n(h) \in \Z$ such that $h=n(h)h_0$. But in that case $d (D, D+h) = n(h) d (D, D+h_0)$, which cannot be $\leq \epsilon $ for all $h \in S_{\epsilon}$, since there are infinitely many elements in $S_{\epsilon}$. This contradiction proves that  the elements of $S'_{\epsilon}$ cannot share the same kernel. Therefore they generate $H^1(\T^2, \R)$. This finishes the proof of Theorem \ref{tore de dim 2}.
\qed
{\small

\bigskip

\noindent

D\'epartement de Math\'ematiques, Universit\'e Montpellier 2, France\\
e-mail : massart@math.univ-montp2.fr
}

\end{document}